\documentclass[preprint,10pt]{elsarticle}
\usepackage[english]{babel}
\usepackage{graphicx}
\usepackage[T1]{fontenc}
\usepackage[utf8x]{inputenc}
\usepackage{color}
\usepackage{tikz}
\usepackage[colorlinks=true]{hyperref}
\usepackage{amsmath}
\usepackage{amssymb}
\usepackage{amsthm}
\usepackage{multirow}
\usepackage{textcomp}
\usepackage{algorithm, algpseudocode}
\usepackage{enumerate}
\usetikzlibrary{positioning}
\usetikzlibrary{decorations.markings}
\usetikzlibrary{decorations.text}
\usetikzlibrary{arrows} 

\usepackage[modulo]{lineno}

\makeatletter
\def\ps@pprintTitle{%
    \let\@oddhead\@empty
    \let\@evenhead\@empty
    \def\@oddfoot{\footnotesize\itshape
         {A preprint} \hfill {October 13, 2022}}%
    \let\@evenfoot\@oddfoot
    }
\makeatother

\begin{document}

\begin{frontmatter}



%\title{Improved Moore-Penrose continuation algorithm for the computation of complex problems}
\title{Improved Moore-Penrose continuation algorithm for the computation of problems with critical points}
% with critical points (of challenging problems)}


\author{S. Léger$^1$, P. Larocque$^1$, D. LeBlanc$^1$}

\address{$^1$D{\'e}partement de math{\'e}matiques
         et de statistique,
         Pavillon Rémi-Rossignol,\\
         18 avenue Antonine-Maillet,
         Universit{\'e} de Moncton,
         Moncton, Canada,
         E1A 3E9.\\ }

\begin{abstract}
Using typical solution strategies to compute the solution curve of challenging problems often leads to the break down of the algorithm. To improve the solution process, numerical continuation methods have proved to be a very efficient tool. However, these methods can still lead to undesired results. In particular, near severe limit points and cusps, the solution process frequently encounters one of the following situations : divergence of the algorithm, a change in direction which makes the algorithm backtrack on a part of the solution curve that has already been obtained and omitting important regions of the solution curve by converging to a point that is much farther than the one anticipated. Detecting these situations is not an easy task when solving practical problems since the shape of the solution curve is not known in advance. This paper will therefore present a modified Moore-Penrose continuation method that will include two key aspects to solve challenging problems : detection of problematic regions during the solution process and additional steps to deal with them. The proposed approach can either be used as a basic continuation method or simply activated when difficulties occur. Numerical examples will be presented to show the efficiency of the new approach. 

\end{abstract}


\end{frontmatter}

\section{INTRODUCTION}
Numerical continuation methods have proved to be a very powerful tool when solving systems of parameterized nonlinear equations. They are widely used in many fields not only to compute solutions that might be hard to obtain otherwise, but also to better understand the physical properties of the problem we are solving. In the case of large deformation problems for example, as the problems are generally driven by a loading parameter, it is not unusual to encounter a section of the solution curve that varies extremely rapidly for some values of the loading parameter due to geometric and/or material non linearities. By using a typical solution strategy, which is often based on Newton like methods, the solution process will usually break down near these difficult regions. Numerical continuation methods do help with these types of simulations; however, some difficulties still occur, particularly near severe limit points and cusps.

Numerical continuation methods can be roughly divided into two main branches, predictor-corrector methods and piecewise linear methods (see \cite{AllGeo1990}). Even if both types of methods share many common features and can be numerically implemented in similar ways, predictor-corrector methods generally perform best when high accuracy is needed. The general idea of these methods is to add a constraint condition to the set of nonlinear equations from which the unknown loading parameter can be determined. As path-following methods are well established, many different variants of these methods are available in the literature, but the ones based on the arc-length method developed by Riks (\cite{Rik1979}) and later modified by Crisfield (\cite{Cri1981}) seem to be the most popular. Dealing with solution curves with critical points, which are commonly categorized into bifurcation points and limit points, is however still a very challenging problem. Many authors have proposed alternative path-following techniques to overcome the problems associated with the convergence issues encountered in the vicinity of critical points (e.g. \cite{Fri1984}, \cite{VanProSim2013}, \cite{PohRamBis2014}). Typical techniques found in the literature are often based on arc-length methods or on methods that add constraint functions based on physical quantities such as total strain (\cite{CheSch1990}) or dissipated energy (\cite{Gut2004}). Even if these techniques improve certain aspects of the basic algorithm, they do not usually fix all problematics and robust algorithms capable of tracing highly nonlinear load-displacement paths are still needed.

As the Moore-Penrose (also known as Gauss-Newton) continuation method, which is a predictor-corrector method, has proved to be very robust when solving challenging problems (\cite{LegForTib2014}, \cite{LegDeiFor2015}, \cite{LegPep2016}), and can easily be implemented in an existing finite element code, this paper will focus on improving this algorithm so that it can better deal with highly nonlinear equilibrium paths.

Similarly to other predictor-corrector methods, the parameter becomes an additional unknown of the problem and can be controlled during the simulation in order to better follow the solution curve, particularly in critical regions where the solution curve is difficult to compute otherwise. In many cases, the solution curve can be computed successfully with the use of the standard Moore-Penrose continuation method. However, this is not always the case. Some typical situations that can be encountered are the divergence of the solution process at critical points, a change in direction at these points which makes the algorithm backtrack on a part of the solution curve that has already been obtained and passing over certain important regions of the curve by converging to a point that is much farther than the one anticipated (which could for example lead to the non tracing of sharp snap-back or snap-through phenomena of structures). The goal of this paper is to improve the Moore-Penrose continuation method so that it performs better in those situations.

To do so, our proposed strategy will include four key ingredients : a recently proposed bifurcation analysis technique called the deflated continuation method (\cite{FarBirFun2015}, \cite{FarBeeBir2016}) to help determine the general behaviour of the solution curve, activation of an angle control during the simulation, a control on the distance between consecutive converged points as well as monitoring a specific component of the tangent vector at the converged points. As the solution curve is generally unknown for practical problems, the goal of the deflated continuation method is to identify distinct solutions for a given parameter value, which will help us detect upcoming difficult regions of the solution curve.  
An additional angle control will also help to deal with the difficult regions. Ligursk\'{y} and Renard (\cite{LigRen2014}) have also improved the Moore-Penrose continuation method by adding an angle control to the algorithm; however, from our tests, we have seen that even if their new algorithm does lead to improved results, it still fails in certain situations. Our proposed strategy therefore includes other features to make it more robust. To show how our proposed strategy is able to address the difficulties observed using other methods, specific validation tests will be considered in this paper. As the initial goal is to better understand the behaviour of the standard algorithm during the simulation, we will consider validation tests for which an analytical solution is known. These tests will incorporate different types of critical points and reproduce solution curves which could be encountered when solving complex practical finite element problems.

The paper is organized as follows. Section~\ref{sec:cont} describes the standard implementation technique for the Moore-Penrose continuation method. Section~\ref{sec:numex} shows specific examples where the standard Moore-Penrose continuation method fails. Section~\ref{sec:keying} presents the key ingredients used in our strategy while Section~\ref{sec:mpimproved} describes the improved algorithm for the Moore-Penrose continuation method. Finally, Section~\ref{sec:validation} is devoted to the validation of this improved strategy.

\section{MOORE-PENROSE CONTINUATION METHOD}\label{sec:cont}
Let us consider a smooth function $F : \mathbb{R}^{N+1} \rightarrow \mathbb{R}^{N}$ for which we want to compute a solution curve of the system
$$F(x)=0$$
starting from a given point on the solution curve. As the system generally depends on one parameter, continuation methods explicitly introduce this parameter in the system of equations as follows :
$$F(x)=F(u, \lambda)=0$$

\noindent The vector of unknowns, $x$, therefore consists of the unknowns of the system, $u$, plus the parameter which is denoted as $\lambda$. To give an example, in the case of large deformation problems, $\lambda$ represents the loading parameter corresponding to either external forces or prescribed displacements or both.

The Moore-Penrose continuation method is a predictor-corrector method and can be summarized as follows (see \cite{DhoGovKuz2003}, \cite{LegDeiFor2015}). Starting from a known point $x^{(i)} \in \mathbb{R}^{N+1}$ on the solution curve, and given a tangent vector $v^{(i)}$ at that point, the next point $x^{(i+1)} \in \mathbb{R}^{N+1}$ on the solution curve as well as its tangent vector $v^{(i+1)}$ can be obtained using the following algorithm :

\begin{itemize}
        \item[\textbullet] $X^0=x^{(i)} + hv^{(i)}$ 
       
       \item[\textbullet] $V^0 = v^{(i)}$
 
 \item[\textbullet] For $k=0,1,2,\ldots,k_{max}$
    \begin{enumerate}
    \item Solve the linear system:
    \begin{equation} \label{sys:step1}
      \begin{cases}
      A(X^k) \, \delta_x^k & = F(X^k)\\
      {V^k}^{\top} \delta_x^k & = 0
      \end{cases}
    \end{equation}
    
    \item Solve the linear system: 
    \begin{equation} \label{sys:step2}
      \begin{cases}
      A(X^k) T^k & = A(X^k)V^k\\
      {V^k}^{\top} T^k & = 0
      \end{cases}
    \end{equation}

    \item Update the solution vector:
    $$X^{k+1}=X^k-\delta_x^k$$
    
    \item Update the tangent vector: 
    $$Z=V^k-T^k$$
    
    \item Normalization of the tangent vector:
    $$V^{k+1}=\frac{Z}{\|Z\|}$$
    
    \item If $\|F(X^k)\|\leq{\varepsilon}_F$ and $\|X^{k+1}-X^k\|\leq{\varepsilon}_x$, convergence attained: 
     \begin{equation*}
      x^{(i+1)}=X^{k+1}, \quad v^{(i+1)}=V^{k+1}
     \end{equation*}

   \end{enumerate} 

\end{itemize}

Let us note that $X^0$ is the prediction point obtained by making a prediction step of length $h$ in the tangential direction, $A(X^k)=[F^{\prime}_u(X^k) \,\,\, F^{\prime}_{\lambda}(X^k)]$ is a rectangular matrix of dimension $n \times (n+1)$, 
$\delta_x^k=[\delta_u^k \,\,\, \delta_{\lambda}^k]^{\top}$ is a correction vector of dimension $(n+1) \times 1$, $k_{max}$ represents the maximum number of iterations allowed while ${\varepsilon}_F$ and ${\varepsilon}_x$ are the desired tolerances on $F$ and $x$ respectively. As for $T^k$, it represents a correction on the tangent vector. 

To help with the convergence of the algorithm, a criterion based on the number of iterations needed for convergence is used to determine if the value of $h$ should be modified or not, and if so, in which way (i.e. increased or decreased). If the algorithm converges to the next point on the solution curve in less than $K_{\textrm{min}}$ iterations, then the step size is increased by a factor of $h_{\textrm{inc}}$. On the other hand, if the algorithm needs more than $K_{\textrm{max}}$ iterations to converge to the next point, the step size is decreased by a factor of $h_{\textrm{dec}}$. Evidently, the chosen values for $K_{\textrm{min}}$ and $K_{\textrm{max}}$ need to be smaller than $k_{\textrm{max}}$ and should be such that $K_{\textrm{min}} < K_{\textrm{max}}$. The smallest value that $h$ can take is noted as $h_{\textrm{min}}$ and once it is attained, if the algorithm is not able to converge to the next point on the solution curve, the remaining of the simulation cannot be completed using the standard Moore-Penrose continuation algorithm. More details can be found in~\cite{LegDeiFor2015}.

\section{NUMERICAL EXAMPLES}\label{sec:numex}
Continuation methods are key to solving difficult problems. However, these types of methods can still lead to unexpected behaviours during the simulation. The goal of this section will be to present some of the situations that can be encountered when using the standard Moore-Penrose continuation method. 

To better analyse the results, the first examples will be based on problems for which the solutions can be obtained analytically. The goal will be to reproduce extreme situations that can occur in more challenging problems using simpler system of equations to better illustrate the behaviours observed. If the continuation method fails in these cases, it will certainly fail when similar curve behaviours are encountered in larger systems of equations. 

In Sections~\ref{subsec:case1} to \ref{subsec:finelemprob}, the following parameter values will be used in the standard Moore-Penrose continuation algorithm : $k_{\textrm{max}}=20$, $K_{\textrm{min}}=5$, $K_{\textrm{max}}=10$, $h_{\textrm{dec}}=0.5$, $h_{\textrm{inc}}=1.5$, $h_{\textrm{min}}=1 \times 10^{-4}$, ${\varepsilon}_F = 1 \times 10^{-7}$ and ${\varepsilon}_x = 1 \times 10^{-7}$.

\subsection{Case I : Divergence of the solution process at critical points}\label{subsec:case1}
At critical points, it is not unusual for the solution process to break down. In these cases, reducing the length of the prediction step, $h$, does not help and the next point on the solution curve cannot be obtained. 

\begin{figure}[!htbp]
\begin{center}
  \begin{tabular}{cc}
   \includegraphics[height=5cm]{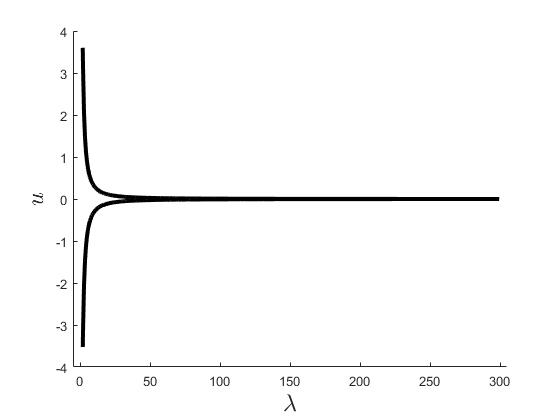} &
   \includegraphics[height=5cm]{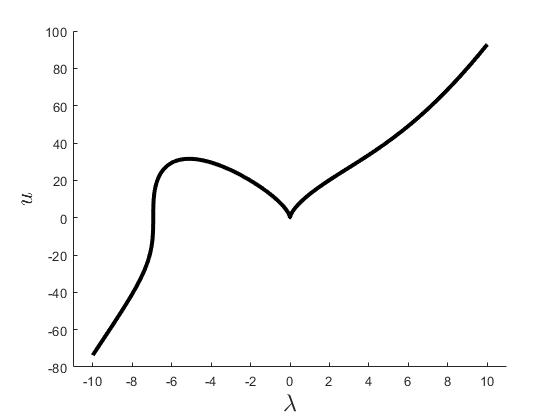} \\
   a) Function $F_a$ & b) Function $F_b$ 
   
  \end{tabular}
  \end{center} 
\caption{Solution curves for functions $F_a$ and $F_b$}
\label{fig:SolutionExacte_fonctionFaFb}           
\end{figure}

To illustrate this behaviour, let us consider two functions $F_a, F_b :  \mathbb{R}^2 \longrightarrow \mathbb{R}$ defined respectively by 
$$F_a(u, \lambda)=-u^2 \lambda^3 - \frac{\lambda}{3}+100 \quad \quad ; \quad \quad F_b(u, \lambda)=2000 \lambda^2 - u^3 + 6 \lambda^5 $$
The solution curves for $F_a(u, \lambda)=0$ and $F_b(u, \lambda)=0$ represent respectively an example with a horizontal limit point and a vertical cusp. Figure~\ref{fig:SolutionExacte_fonctionFaFb} illustrates the solution curve for both examples.  

\begin{figure}[!htbp]
\begin{center}
   \includegraphics[height=6cm]{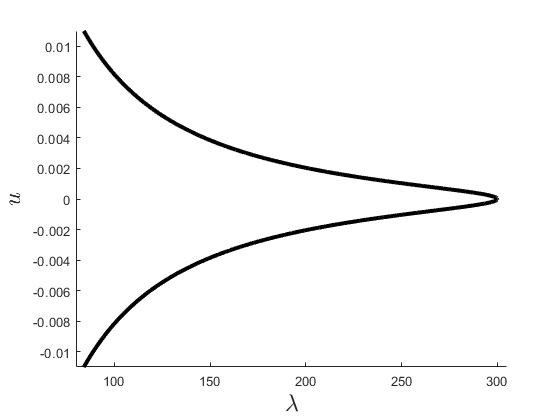}
    \end{center} 
\caption{Solution curve for $F_a(u, \lambda)=0$ : zoom in the region of the horizontal limit point}
\label{fig:SolutionExacteZoom_fonctionFa}           
\end{figure}

As can be seen in Figure~\ref{fig:SolutionExacte_fonctionFaFb}a), the horizontal limit point is very severe. This is done to reproduce situations that are encountered in challenging applications. Figure~\ref{fig:SolutionExacteZoom_fonctionFa} illustrates a zoom in the critical region.

\begin{figure}[!htbp]
\begin{center}
   \includegraphics[height=5cm]{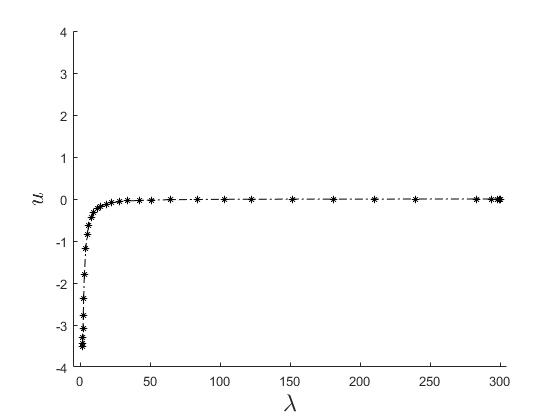} \end{center} 
\caption{Numerical solution obtained in the case of $F_a$ using the standard Moore-Penrose continuation method}
\label{fig:SolutionNum_fonctionFa}           
\end{figure}

\begin{figure}[!htbp]
\begin{center}
   \includegraphics[height=6cm]{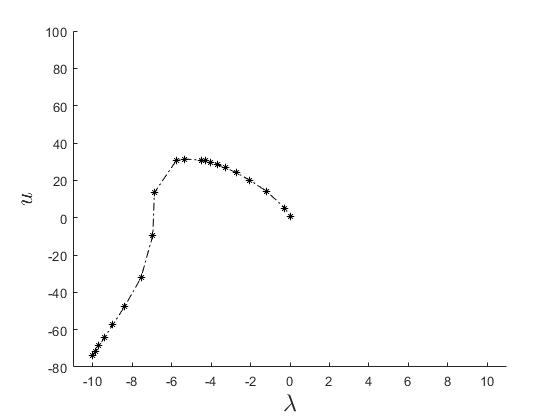}
    \end{center} 
\caption{Numerical solution obtained in the case of $F_b$ using the standard Moore-Penrose continuation method}
\label{fig:SolutionNum_fonctionFb}           
\end{figure}

By applying the standard Moore-Penrose continuation method, the solution process breaks down as can be seen in Figures~\ref{fig:SolutionNum_fonctionFa} and \ref{fig:SolutionNum_fonctionFb}. In both cases, the algorithm fails at the critical point and is unable to pursue with the simulation.

\subsection{Case II : Convergence towards a part of the solution curve that has already been computed}\label{subsec:case2}
Another issue that can occur at critical points is that when the algorithm has difficulty to converge to the next point on the solution curve it ends up converging to a point on a region of the curve that has already been computed. Since the tangent vector now points in the opposite direction, the algorithm then proceeds to backtrack on the solution curve, reproducing the same part of the solution curve that has just been computed. One of the major difficulties of this situation is detecting numerically that the algorithm has changed direction since a change in the direction of the tangent vector does not automatically mean that the algorithm is in backtracking mode.

\begin{figure}[!htbp]
\begin{center}
  \begin{tabular}{cc}
   \includegraphics[height=5cm]{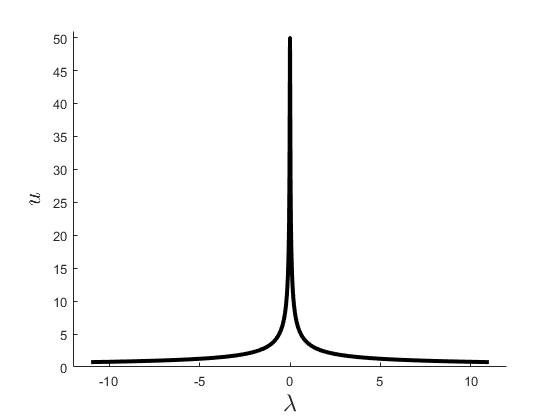} &
   \includegraphics[height=5cm]{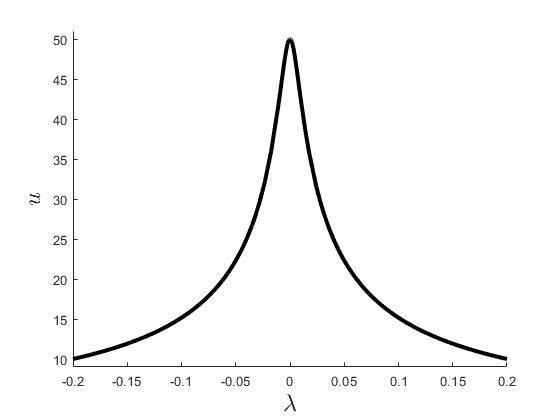} \\
   a) Overall view & b) Zoom in the critical region 
   
  \end{tabular}
  \end{center} 
\caption{Solution curve for function $F_c$}
\label{fig:SolutionExacte_fonctionFc}           
\end{figure}

To illustrate this behaviour, let us consider a function $F_c :  \mathbb{R}^2 \longrightarrow \mathbb{R}$ defined by
$$F_c(u,\lambda)=-u^3 \lambda^2 - u + 50$$
for which the solution curve of $F_c(u,\lambda)=0$ represents an example with a vertical limit point. The solution curve for this equation is shown in Figure~\ref{fig:SolutionExacte_fonctionFc}.

\begin{figure}[!htbp]
\begin{center}
  \begin{tabular}{cc}
   \includegraphics[height=5cm]{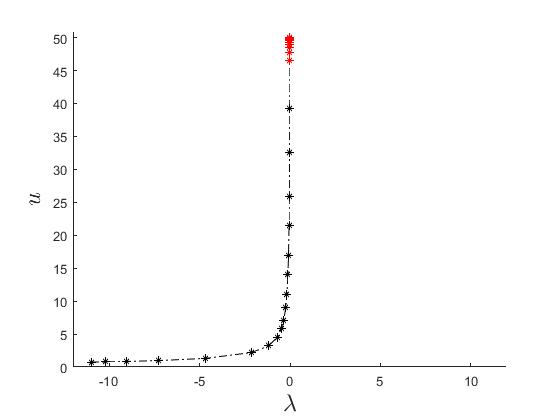} &
   \includegraphics[height=5cm]{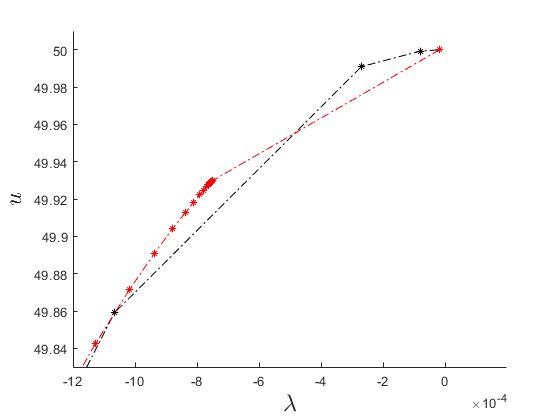} \\
   a) Overall view & b) Zoom in the critical region 
   
  \end{tabular}
  \end{center} 
\caption{Numerical solution obtained in the case of $F_c$ using the standard Moore-Penrose continuation method}
\label{fig:SolutionNum_fonctionFc}           
\end{figure}

As can be seen in Figure~\ref{fig:SolutionNum_fonctionFc}, when the standard Moore-Penrose continuation method reaches the critical point, it starts to backtrack on the initial part of the solution curve. To illustrate this behaviour, the converged points obtained after the algorithm reaches the critical point are shown in red. By comparing the black converged points with the red ones in Figure ~\ref{fig:SolutionNum_fonctionFc}b), it might seem that the algorithm is not landing exactly on the same part of the curve. However, this is misleading as linear segments are drawn between converged points. All points (red and black) are in fact located on the same part of the solution curve. 

In this specific case, by modifying certain values of the parameters (e.g. taking $h_{\textrm{min}}$ extremely small), we were able to complete the simulation. However, it was extremely hard to find parameters for which the simulation was successful. This is far from being an ideal situation since very robust numerical algorithms are needed when solving practical problems.

\subsection{Case III : Convergence to a point on the solution curve that is much farther than the one anticipated}\label{subsec:case3}
Converging to a point far along the solution curve is not automatically an issue, unless it means that an important region of the curve has been bypassed. Limit points, cusps and bifurcation points are hard to deal with numerically, but knowing that some of these points exist for a particular problem can be important as it helps us better understand the physical properties of the problem we are solving. 

\begin{figure}[!htbp]
\begin{center}
  \begin{tabular}{c}
   \includegraphics[height=5cm]{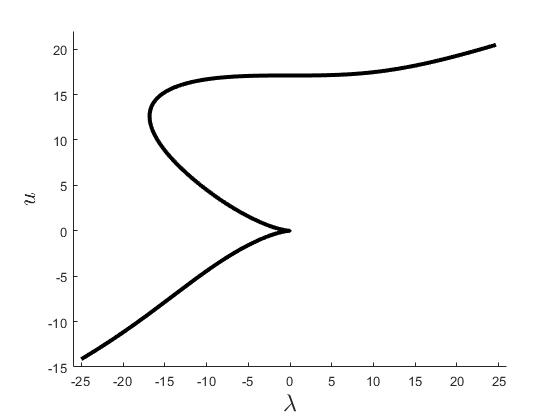} 
     \end{tabular}
  \end{center} 
\caption{Solution curve for function $F_d$}
\label{fig:SolutionExacte_fonctionFd}           
\end{figure}

To illustrate this behaviour, let us consider a function $F_d :  \mathbb{R}^2 \longrightarrow \mathbb{R}$ defined by
$$F_d(u,\lambda)=-500u^2 - 10\lambda^3 + \frac{1}{10}u^5$$
for which $F_d(u,\lambda)=0$ represents an example with a horizontal cusp. Figure~\ref{fig:SolutionExacte_fonctionFd} illustrates the solution curve for this example.

\begin{figure}[!htbp]
\begin{center}
  \begin{tabular}{c}
   \includegraphics[height=5cm]{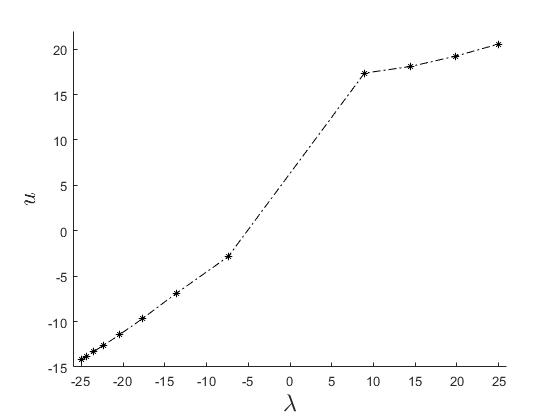} 
     \end{tabular}
  \end{center} 
\caption{Numerical solution obtained in the case of $F_d$ using the standard Moore-Penrose continuation method}
\label{fig:SolutionNum_fonctionFd}       
\end{figure}

As can be seen in Figure~\ref{fig:SolutionNum_fonctionFd}, the standard Moore-Penrose continuation method neglects the entire difficult region (i.e. horizontal cusp) of the solution curve by converging to a point that is much farther along the path.

\subsection{Finite element problems}\label{subsec:finelemprob}
\subsubsection{Bratu problem}
The difficulties encountered in Cases I to III frequently occur when solving certain types of finite element problems. To give an example, let us consider the classical Bratu problem~\cite{Bra1914}. As this problem appears in a large variety of application areas (see~\cite{JacSch2002} for a summary of the history of the problem), it is commonly used as a test problem for numerical methods. This problem, which is an elliptic nonlinear partial differential equation with homogeneous Dirichlet boundary conditions, is given by

$$
\begin{cases}
\Delta w + \lambda e^{w} = 0 \quad \textrm{in} \quad \Omega \\
w=0 \quad \textrm{on} \quad \Gamma
\end{cases}
$$

\noindent where $\lambda>0$. $\Omega$ represents the bounded domain with boundary $\Gamma$. In the one dimensional case, the problem reduces to

$$
\begin{cases}
w_{xx} + \lambda e^{w} = 0, \quad 0 \le x \le 1 \\
w(0)=w(1)=0
\end{cases}
$$

\begin{figure}[!htbp]
\begin{center}
   \includegraphics[height=5cm]{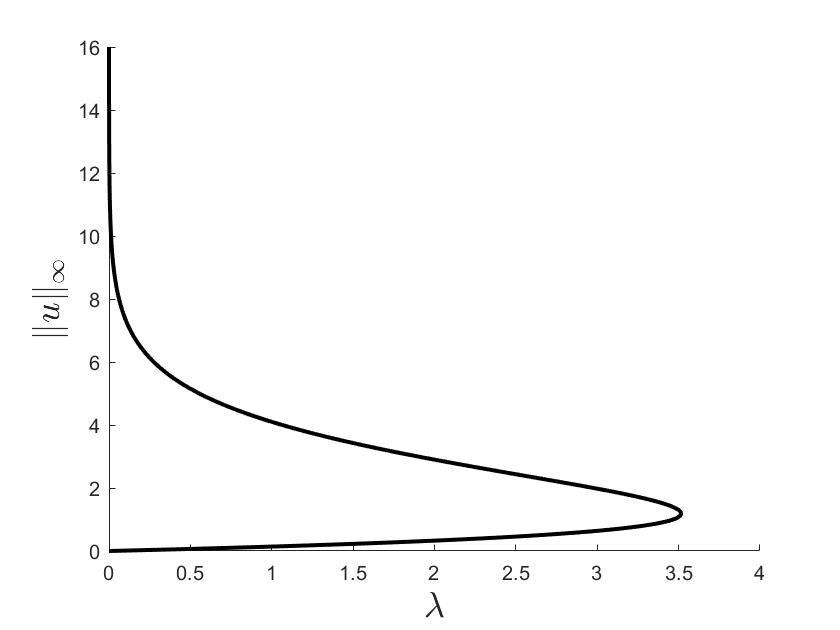} 
 \end{center} 
\caption{Solution curve for the Bratu problem}
\label{fig:SolutionExactBrattu}           
\end{figure}

\noindent The exact solution of this problem is illustrated in Figure~\ref{fig:SolutionExactBrattu}. As can be seen, it exhibits a horizontal limit point. To make this limit point even more severe, we can simply set $w=\gamma u$ with $\gamma$ a constant, which leads to 

$$
\begin{cases}
\gamma u_{xx} + \lambda e^{\gamma u} = 0, \quad 0 \le x \le 1 \\
u(0)=u(1)=0
\end{cases}
$$

\noindent The solution curve of this modified problem is compressed vertically by a factor of $\gamma$.

\begin{figure}[!htbp]
\begin{center}
  \begin{tabular}{cc}
   \includegraphics[height=5cm]{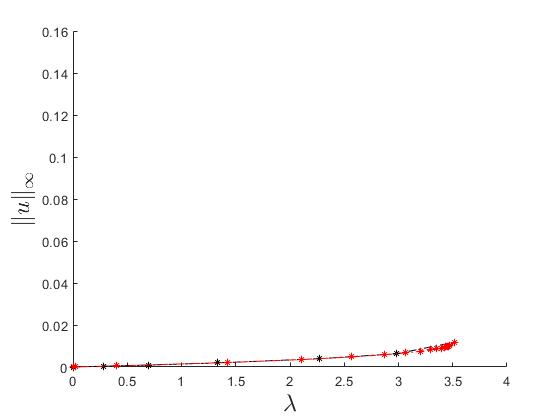} &
   \includegraphics[height=5cm]{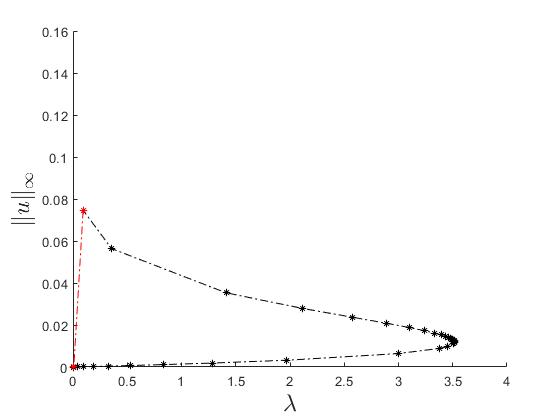} \\
   a) Scenario 1 & b) Scenario 2 
   
  \end{tabular}
  \end{center} 
\caption{Two examples of the numerical solution obtained in the case of the modified Bratu problem for $\gamma=100$ when using the standard Moore-Penrose continuation method}
\label{fig:SolutionNum_BratuK100_MPBase1}           
\end{figure}

When solving the modified Bratu problem with $\gamma=100$ using the finite element method (with quadratic elements) and the standard Moore-Penrose continuation method, it is not unusual for the algorithm to fail, all depending on the choice of the initial point and parameter values. Figure~\ref{fig:SolutionNum_BratuK100_MPBase1} shows two examples of situations that were encountered during the solution process. In the first scenario, once the algorithm reaches the critical region, it jumps to a point on the first part of the solution curve and then begins backtracking. The initial points are shown in black while the converged points obtained after the algorithm has reached the critical points are shown in red. As for the second scenario, the algorithm is able to bypass the extremity of the limit point, but then jumps to the initial point later on in the simulation. This jump is illustrated in red. The choice of a slightly different initial step length $h$ explains the difference in the scenarios.

\subsubsection{Other finite element problem}
Let us now consider the following one dimensional non linear problem :
$$\begin{cases}
(u(x))^{\alpha} - \nabla \cdot (q(x) \nabla u(x)) = r(x) \quad \textrm{in} \quad \Omega \\
u(x)=0 \quad \textrm{on} \quad \Gamma
\end{cases}$$

To test the standard Moore-Penrose continuation method in other finite element applications where difficult regions appear on the solution curve, we can use the method of manufactured solution as described in~\cite{ChaForFor2010}. The idea of this method is simple. It simply consists of injecting an analytical expression for $u(x)$ in the differential equation to generate the function $r(x)$ which is then used as an artificial source term in the finite element code. The numerical solution can then be compared with the analytical one. 

\begin{figure}[!htbp]
\begin{center}
   \includegraphics[height=5cm]{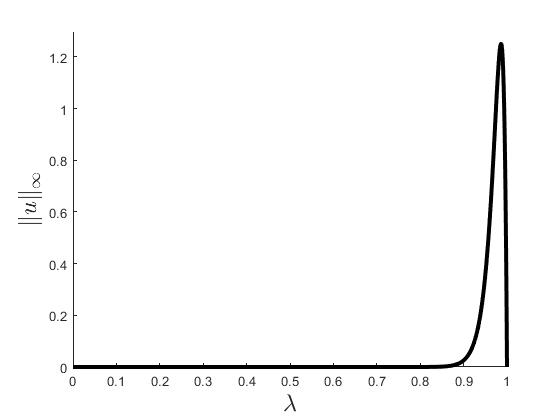} 
 \end{center} 
\caption{Solution curve for the manufactured solution problem}
\label{fig:SolutionExactManuf}           
\end{figure}

Let us consider the function 
$$u(x)=\zeta \lambda^{\eta} (1-\lambda^{\eta})(1-x)x$$
with $\zeta$ and $\eta$ constant for the analytical solution and set $\Omega=[0, 1]$, $\alpha=2$ and $q(x)=1$. This solution exhibits a vertical limit point and its severity is determined by the values of $\zeta$ and $\eta$. The value of $\zeta$ influences the height of the limit point while the value of $\eta$ influences its width. Figure~\ref{fig:SolutionExactManuf} illustrates the solution in the case of $\zeta=20$ and $\eta=50$.

\begin{figure}[!htbp]
\begin{center}
   \includegraphics[height=5cm]{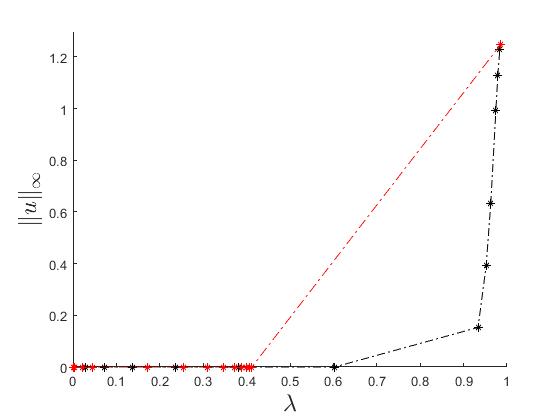} 
 \end{center} 
\caption{Example of the numerical solution obtained in the case of the manufactured solution problem when using the standard Moore-Penrose continuation method}
\label{fig:SolutionAlgoBaseManuf}           
\end{figure}

By solving this finite element problem with the help of the standard Moore-Penrose continuation method, the algorithm again leads to many undesired results depending on the choice of the initial point and parameter values. Figure~\ref{fig:SolutionAlgoBaseManuf} illustrates one of these situations in which case the algorithm starts to backtrack on the solution curve once it reaches the top part of the critical region. The points shown in red are the converged points obtained after the critical point is reached. 

The examples presented in Section~\ref{sec:numex} will serve as a reference for our improved Moore-Penrose continuation algorithm.

\section{KEY INGREDIENTS}\label{sec:keying}
Our improved strategy is based on the addition of four ingredients in the standard Moore-Penrose continuation method : the use of the deflated continuation algorithm (\cite{FarBirFun2015}, \cite{FarBeeBir2016}), an additional angle control (\cite{LigRen2014}), a control on the distance between consecutive converged points and monitoring the sign of $v_{\lambda}$. All of these tools are discussed below.

\subsection{Deflated continuation}
In~\cite{FarBirFun2015}, Farrell and $al.$ present a new algorithm to compute bifurcation diagrams. This new approach, which combines classical continuation with a deflation technique, is called deflated continuation and is able to discover other solution branches from known ones. The key aspect of this algorithm is the deflation technique, which systematically modifies the nonlinear problem so that Newton's method does not converge to a solution that has already been obtained. The unknown solutions are obtained using the same initial guess.

In our case, the goal is not necessarily to compute the bifurcation diagram in its entirety, but rather to know if we could be dealing with critical points during the simulation. Knowing that multiple branches exist (or more generally that multiple solutions exist for a fixed value of $\lambda$) enables us to better understand the physical properties of the problem and to adjust the Moore-Penrose continuation algorithm accordingly in order to avoid numerical difficulties. The deflated continuation algorithm can be described as follows. More details can be found in \cite{FarBirFun2015} and \cite{FarBeeBir2016}.

Starting from a fixed value for the parameter $\lambda$, the goal is to find a solution to the equation
$$F(u,\lambda) = f(u)=0$$
where $f : \mathbb{R}^{N} \rightarrow \mathbb{R}^{N}$ by starting from an initial guess $u_0$. The solution obtained is denoted by $u_1^{\ast}$. To determine if other solutions exist for that value of $\lambda$, the following modified problem is then constructed :
$$G(u)=M(u; u_1^{\ast})f(u)$$
where $M(u; u_1^{\ast})$ is the deflation operator to the residual $f$ constructed in such a way that the deflated residual will satisfy the following two properties :
\begin{itemize}
\item preservation of the solutions of $f$ (i.e. for $u \ne u_1^{\ast}$, $G(u)=0$ if and only if $f(u)=0)$
\item applying Newton's method to $G$ will not find solution $u_1^{\ast}$ again
\end{itemize}

In the work presented by Farrell and $al.$, the shifted deflation operator is used :
$$M(u; u_1^{\ast}) = \left ( \frac{1}{||u-u_1^{\ast}||^p} + \sigma \right) I$$
where $I$ is the identity matrix in $\mathbb{R}^{N}$, $p$ is the power, and $\sigma$ is the shift. In our work, the same deflation operator will be used with $p=2$ and $\sigma=1$.

If another solution is found for the fixed parameter value $\lambda$, the procedure is repeated by deflating once again the residual. For example, if deflation leads to solution $u_2^{\ast}$, the following modified problem will then be constructed :
$$H(u)=M(u; u_2^{\ast})G(u)$$
and the search will continue in order to see if other solutions can be found. When attempting to find a new solution, the solutions that have already been found for the previous $\lambda$ value are used, in turn, as an initial guess for Newton's method until a new solution is found. If no other solutions can be found, the parameter $\lambda$ is increased to its next value.

In our algorithm, the deflation technique is not applied at every iteration, but rather at every $N$ iterations of the Moore-Penrose continuation method to determine if other solution branches exist for that value of $\lambda$. The knowledge obtained from these steps enables us to determine what approach should be used if difficulties are encountered during the simulation. To reduce the computational cost, we do not attempt to find all possible branches at all times. The first time that an additional solution is found for a fixed value of $\lambda$, we wait until the next deflated step to make more attempts at finding solutions. At this next step, both solutions that were found previously are used, in turn, to find solutions for the current value of $\lambda$. The same procedure is applied for the other deflation steps : all solutions that were found at the last deflation step are used in turn to find the solutions at the current step.

\subsection{Angle control}
The idea of adding an angle control to our algorithm came after reading the work by Ligursk\'{y} and Renard (\cite{LigRen2014}) in which they describe a piecewise-smooth inexact Moore-Penrose predictor-corrector algorithm. In this algorithm, there is a parameter $c_{\textrm{min}}$ which controls changes of direction between tangent vectors at two consecutive converged points. If the deviation between the newly computed tangent and the previous one is too large or if the corrector step does not converge, the step length $h$ used in the prediction step is reduced and another attempt is made by starting from the last converged point. On the other hand, the step length can also be increased if certain conditions are respected. 

As was done in \cite{LigRen2014}, if the angle criteria is not respected, that is if 
$$(v^{(i+1)})^{\top} v^{(i)} <  c_{\textrm{min}}$$
the step length is reduced as follows :
$$h := \max\{h_{\textrm{dec}} h, h_{\textrm{min}}\}$$
where $h_{\textrm{dec}}$ is the scale factor for the reduction of $h$. Evidently, $h$ cannot be reduced below $h_{\textrm{min}}$.

\subsection{Distance between consecutive converged points}\label{subsec:distance}
In order to make sure that important sections of the solution curve are not bypassed, we have also included another ingredient in our algorithm, which consists in controlling the distance between consecutive converged points ($x^{(i)}$ and $x^{(i+1)}$), both in terms of $u$ and $\lambda$. The distance between these two points is monitored throughout the simulation. At a certain point, if the distance in terms of $||u||$ is greater than a fixed parameter value $\delta_{\textrm{maxU}}$ or if the distance in terms of $\lambda$ is greater than a fixed parameter value $\delta_{\textrm{maxL}}$, the jump is considered to be too great and another attempt to find $x^{(i+1)}$ is made by returning to the point $x^{(i)}$ and reducing the step length $h$ by a factor of $h_{\textrm{dec}}$. If the point obtained using this reduced step length is still too far along the solution curve, the process is repeated until $h$ reaches $h_{\textrm{min}}$. In the case where a point $x^{(i+1)}$ that respects the desired distances cannot be found for $h \ge h_{\textrm{min}}$, a special procedure will be applied. More details will be provided in Section~\ref{sec:mpimproved}.

\subsection{Sign of $v_{\lambda}$}\label{subsec:sign_v_lambda}
A change in the sign of $v_{\lambda}$ between two consecutive converged points could mean, amongst other scenarios, that we are dealing with a horizontal limit point or cusp or that the algorithm is starting to backtrack on the solution curve. The sign of $v_{\lambda}$ therefore needs to be monitored throughout the simulation as both of these situations need to be dealt with carefully. The information obtained from the deflation steps will be helpful in determining which steps should be completed to pursue the simulation.

\section{IMPROVED MOORE-PENROSE CONTINUATION METHOD} \label{sec:mpimproved}

Our algorithm is based on the standard Moore-Penrose continuation method, but modifications were made to deal with the numerical difficulties that can occur when solving challenging problems. This section will explain how the key ingredients presented in Section~\ref{sec:keying} are incorporated in our strategy.

The goal of the following algorithm is to be able to deal with the difficulties encountered in Section~\ref{sec:numex} which were caused by severe limit points or cusps (vertical or horizontal). We will distinguish two different cases based on the results obtained.

\subsection{Algorithm}
When difficulties occur during the simulation, they are often due to limit points or cusps. To determine if we are dealing with a horizontal limit point or cusp versus a vertical limit point or cusp, we will use the information obtained during the deflated steps. 

\begin{itemize}
\item \textbf{Vertical limit points or cusps}

In the case of severe vertical limit points or cusps, the deflation technique should, in most cases, only detect one branch for each value of $\lambda$ in its neighbourhood. Other possibilities could include detecting multiple branches but very far from the current one or detecting close branches but for which the distance between them becomes greater when $\lambda$ is increased. In either of these cases, the following additional steps should be applied to the standard Moore-Penrose continuation method once $x^{(i+1)}$ has been computed :

\begin{enumerate}
\item[a)] Verification of the distance, in terms of $u$ and $\lambda$, the sign of $v_{\lambda}$ and the angle between the current and previous converged points. In this case, if $v_{\lambda}^{(i)}$ and $v_{\lambda}^{(i+1)}$ are of opposite signs, it means that the algorithm is starting to backtrack on the solution curve since the deflation technique has ruled out the possibility of dealing with a horizontal limit point or cusp.
\begin{itemize}
\item[\textbullet] While $||u^{(i+1)}-u^{(i)}||>\delta_{\textrm{maxU}}$ or $|\lambda^{(i+1)}-\lambda^{(i)}|>\delta_{\textrm{maxL}}$ or $v_{\lambda}^{(i+1)} v_{\lambda}^{(i)} < 0$ or $(v^{(i+1)})^{\top} v^{(i)} <  c_{\textrm{min}}$, set $h := \max\{h_{\textrm{dec}} h, h_{\textrm{min}}\}$ and return to the last converged point $x^{(i)}$ (with its tangent vector $v^{(i)}$) to make another attempt to find $x^{(i+1)}$. If a point $x^{(i+1)}$ that respects the desired criteria is found for $h \ge h_{\textrm{min}}$, then the algorithm continues normally until other difficulties occur. If not, proceed to step b).
\end{itemize} 

\item[b)] If $h=h_{\textrm{min}}$ and the algorithm is not able to converge to an appropriate point $x^{(i+1)}$, apply the \textit{Vertical turning point method} (see Algorithm~\ref{algo:VertTurnPoint}). The idea of this method is illustrated in Figure~\ref{fig:VertTurnMethod}.
\end{enumerate}

\begin{algorithm}
\caption{Vertical turning point method} \label{algo:VertTurnPoint}
\begin{algorithmic}
\Require $x^{(i)}$, $\Delta_{\lambda}$, $\epsilon_{\lambda}$
\Statex
\Statex Step 1 : Starting from the point $x^{(i)}$, find the solution corresponding to $\lambda=\lambda + \Delta_{\lambda}$ (or $\lambda = \lambda - \Delta_{\lambda}$ if $v_{\lambda}^{(i)}<0$) by using the basic Newton method, where $\Delta_{\lambda}$ represents a small variation of $\lambda$. This new point is denoted as $Z^{\ast}$.
\Statex
\If{$v_{\lambda}$>0}
\State $Z^0=[x_u^{(i)} \,\,\, x_{\lambda}^{(i)}+\Delta_{\lambda}]$
\State $\epsilon_{\lambda}^{\ast} = \epsilon_{\lambda}$
\Else[$v_{\lambda}$<0]
\State $Z^0=[x_u^{(i)} \,\,\, x_{\lambda}^{(i)}-\Delta_{\lambda}]$
\State $\epsilon_{\lambda}^{\ast} = -\epsilon_{\lambda}$
\EndIf
\Statex Using $Z^0$, apply the basic Newton method to find $Z^{\ast}$ such that $F(Z^{\ast})=0$
\Statex
\Statex Step 2 : Approximate the tangent vector at $Z^{\ast}$ by the secant passing between the points $x^{(i)}$ and $Z^{\ast}$ and then normalize this approximated tangent vector. The resulting vector is denoted as $W$.
\Statex
\Statex $W=Z^{\ast}-x^{(i)}$
\Statex $W = W/||W||$
\Statex
\Statex Step 3 : To ensure that the approximation of the new tangent vector is not vertical, add $\epsilon_{\lambda}^{\ast}$ to the $\lambda$ component of $W$ and renormalize the vector to obtain $W^{\ast}$. This is simply done to ensure that the algorithm will continue to pursue the curve in the right direction.
\Statex
\Statex $W^{\ast} = [W_u \,\,\, W_{\lambda}+\epsilon_{\lambda}^{\ast}]$
\Statex $W^{\ast} = W^{\ast}/||W^{\ast}||$
\Statex
\Statex Step 4 : Find the next point on the solution curve by setting $x^{(i)}=Z^{\ast}$ and $v^{(i)}=W^{\ast}$ and applying the standard Moore-Penrose continuation algorithm. The angle control is then deactivated until it needs to be activated again.
 
\end{algorithmic}
\end{algorithm}

\begin{figure}[!htbp]
    \centering
     \begin{tabular}{ccc}
\begin{tikzpicture}[scale = 1.5]
\draw (-2,5) to [out=300, in=95]
node[anchor=east,pos=0.4] {$x^{(i)}$} 
node[fill,circle,inner sep=0pt,minimum size=4pt,pos = 0.4](point){}
node[blue,anchor=east,pos=0.6] {$Z^*$} 
node[blue,fill,circle,inner sep=0pt,minimum size=4pt,pos = 0.6](pointZ){}
node[anchor= east,pos=0.75] {$x^{(i+1)}$} 
node[fill,circle,inner sep=0pt,minimum size=4pt,pos = 0.75](point2){}
node[pos=0.9] (vect2) {}
(0,0);

\draw (2,5) to [out=240, in=85] (0,0);

\path (point)
          edge[-stealth',black] node[anchor=east, left = 1pt] {${{v}^{(i)}}$} (-0.45,2.25); 
\path (point)
          edge[-stealth',blue] node[anchor=east, pos = 0.75] {${W}$} (pointZ); 
\path (point)
          edge[-stealth',red] node[anchor = west] {${{W}^{*}}$} (-0.25,2.25); 
\path (point2)
          edge[-stealth',black] node[anchor = east] {${{v}^{(i+1)}}$} (vect2); 
\end{tikzpicture} &  &
\begin{tikzpicture}[scale = 1.5]
\draw (-2,5) to [out=300, in=95]
node[anchor=east,pos=0.95] {$x^{(i)}$} 
node[fill,circle,inner sep=0pt,minimum size=4pt,pos = 0.95](point){}
(0,0);
\draw (2,5) to [out=240, in=85] 
node[blue,anchor=west,pos=0.7] {$Z^*$} 
node[blue,fill,circle,inner sep=0pt,minimum size=4pt,pos = 0.7](pointZ){}
node[anchor = west,pos=0.65] {$x^{(i+1)}$} 
node[fill,circle,inner sep=0pt,minimum size=4pt,pos = 0.65](point2){}
node[pos=0.45, left = 2pt] (vect2) {}
(0,0);

\path (point)
          edge[-stealth',blue] node[anchor=east] {${W}$} (pointZ); 
\path (point)
          edge[-stealth',black] node[anchor=east] {${{v}^{(i)}}$} (0.1,-0.5); 
\path (point)
          edge[-stealth',red] node[anchor = west] {${{W}^{*}}$} (0.5,1.5); 
\path (point2)
          edge[-stealth',black] node[anchor = east] {${{v}^{(i+1)}}$} (vect2); 
\end{tikzpicture} \\
\end{tabular}
\caption{Illustration of the vertical turning point method}
\label{fig:VertTurnMethod} 
\end{figure}
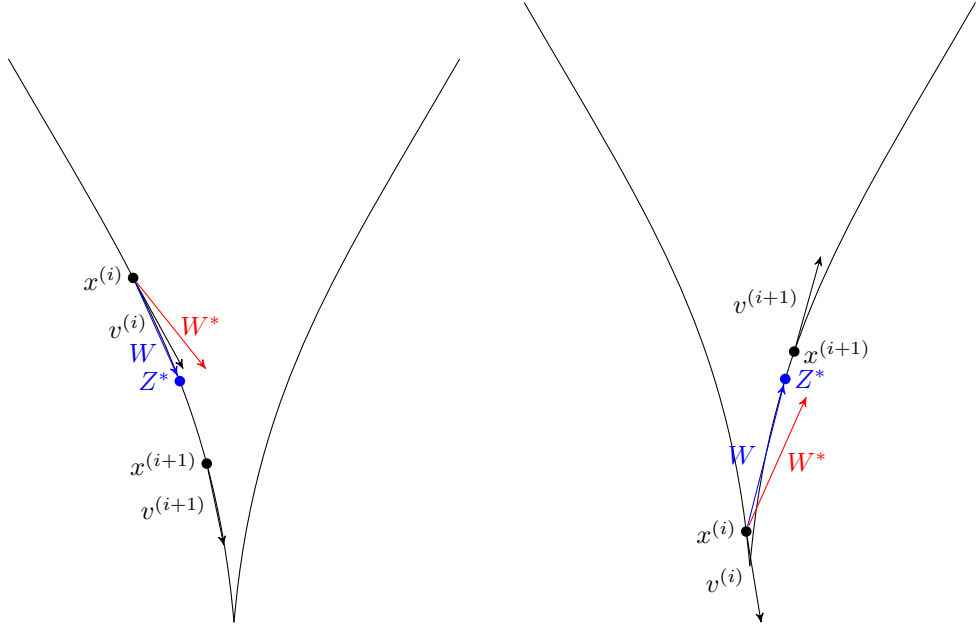

If the standard Moore-Penrose continuation algorithm simply has difficulty converging to the next point on the solution curve (i.e. the algorithm is blocked at point $x^{(i)}$ and reducing the step size does not help), the vertical turning point method would also be applied in the case where deflation has only detected one solution branch (or as discussed earlier when multiple branches are detected but they are far from the current one or when close branches are detected but the distance between them increases with $\lambda$).

\item \textbf{Horizontal limit points or cusps}

If horizontal limit points or cusps are present in the solution curve, multiple branches (or solutions) should be detected, for a fixed value of $\lambda$, by the deflation technique when approaching these critical points. As the numerical difficulties usually occur when the critical points are severe, if multiple branches are detected for a certain value of $\lambda$, our algorithm therefore also monitors the distance between the principal branch and its closest branch. This distance will be noted as $\delta$. As soon as this distance reaches a critical value set by the user, noted as $\delta_{\textrm{crit}}$ in this work, it means that the branches are getting close and that numerical problems could be encountered.

\noindent In the case of multiple branches, there is therefore two possibilities, the first one being the situation where the distance between the branches should not be a cause for concern (i.e. $\delta > \delta_{\textrm{crit}}$ or $\delta < \delta_{\textrm{crit}}$ with $\delta$ increasing with the Moore-Penrose iterations) which has already been treated above and the second one the situation where the branches are getting closer as the simulation continues (i.e. $\delta < \delta_{\textrm{crit}}$ with $\delta$ decreasing with the Moore-Penrose iterations). In this second case, our strategy is to apply the \textit{Horizontal turning point method} (see Algorithm~\ref{algo:HorizTurnPoint}). The idea of this method is illustrated in Figure~\ref{fig:HorizTurnMethod}.

\begin{algorithm}
\caption{Horizontal turning point method} \label{algo:HorizTurnPoint}
\begin{algorithmic} 
\Require $x^{(i)}$, $v^{(i)}$, $Y^{(i)}$, $\delta_{\textrm{crit}}$, $h$, $\epsilon_{\textrm{diff}}$, $\delta_{\textrm{maxU}}$, $\delta_{\textrm{maxL}}$, $c_{\textrm{min}}$, $h_{\textrm{dec}}$, $h_{\textrm{min}}$
\Statex
\Statex Step 1 : When a point $Y^{(i)}$ is detected on the secondary branch at a distance inferior to $\delta_{\textrm{crit}}$ from the current point $x^{(i)}$, approximate the tangent vector at this point.
\Statex
\Statex Find $w^{(i)}$ such that $F'(Y^{(i)})w^{(i)}=0$ and $||w^{(i)}||=1$.
\Statex
\Statex Step 2 : Apply a modified version of the Moore-Penrose continuation method on each branch until both branches coincide or that the computation of a new point on the branch fails. To better describe this strategy, the new points computed on the principal branch will be denoted normally as $x^{(i+1)}$ with tangent vector $v^{(i+1)}$ while those on the secondary branch will be denoted as $Y^{(i+1)}$ with tangent vector $w^{(i+1)}$.
\Statex
\Statex $x_{\textrm{cur}}=x^{(i)}$
\Statex $Y_{\textrm{cur}}=Y^{(i)}$
\Statex
\Statex $h_{\textrm{ConvPrinc}}=h$
\Statex $h_{\textrm{ConvSec}}=h$
\Statex
\Statex $\textrm{ConvPrinc}=\textrm{true}$
\Statex $\textrm{ConvSec}=\textrm{true}$
\Statex
\While{$||x_{\textrm{cur}} - Y_{\textrm{cur}}||>\epsilon_{\textrm{diff}}$ AND ($\textrm{ConvPrinc}=\textrm{true}$ OR $\textrm{ConvSec}=\textrm{true}$)}
\If{$\textrm{ConvPrinc}=\textrm{true}$}
	\State Starting from point $x^{(i)}$, calculate $x^{(i+1)}$ by using the standard 
	\State Moore-Penrose continuation algorithm.
	\Statex
\While{$||x_u^{(i+1)}-x_u^{(i)}||>\delta_{\textrm{maxU}}$ OR $|x_{\lambda}^{(i+1)}-x_{\lambda}^{(i)}|>\delta_{\textrm{maxL}}$ OR $(v^{(i+1)})^{\top} v^{(i)} <  c_{\textrm{min}}$}
\State Set $h_{\textrm{ConvPrinc}} := \max\{h_{\textrm{dec}} h_{\textrm{ConvPrinc}}, h_{\textrm{min}}\}$ and return to the last 
\State converged point $x^{(i)}$ (with its tangent vector $v^{(i)}$) to make another 
\State attempt to find $x^{(i+1)}$. 
\Statex
\If{Appropriate $x^{(i+1)}$ is found for $h_{\textrm{ConvPrinc}} \ge h_{\textrm{min}}$}
\If{$v_{\lambda}^{(i+1)} v_{\lambda}^{(i)} > 0$}
\State $x_{\textrm{cur}}=x^{(i+1)}$
\Else
\State Continuation of principal branch terminated.
\State ConvPrinc=false
\State $x_{\textrm{cur}}=x^{(i)}$
\State $x^{\textrm{final}}=x^{(i)}$
\EndIf
\Else
\State Continuation of principal branch terminated.
\State ConvPrinc=false
\State $x_{\textrm{cur}}=x^{(i)}$
\State $x^{\textrm{final}}=x^{(i)}$
\EndIf
\EndWhile
\EndIf
\algstore{Alg2}
\end{algorithmic}
\end{algorithm}

\begin{algorithm}
\begin{algorithmic} 

\algrestore{Alg2}

\If{$\textrm{ConvSec}=\textrm{true}$}
	\State Starting from point $Y^{(i)}$, calculate $Y^{(i+1)}$ by using the standard 
	\State Moore-Penrose continuation algorithm.
	\Statex
\While{$||Y_u^{(i+1)}-Y_u^{(i)}||>\delta_{\textrm{maxU}}$ OR $|Y_{\lambda}^{(i+1)}-Y_{\lambda}^{(i)}|>\delta_{\textrm{maxL}}$ OR $(w^{(i+1)})^{\top} w^{(i)} <  c_{\textrm{min}}$}
\State Set $h_{\textrm{ConvSec}} := \max\{h_{\textrm{dec}} h_{\textrm{ConvSec}}, h_{\textrm{min}}\}$ and return to the last 
 \State converged point $Y^{(i)}$ (with its tangent vector $w^{(i)}$) to make another  
\State attempt to find $Y^{(i+1)}$. 
\Statex
\If{Appropriate $Y^{(i+1)}$ is found for $h_{\textrm{ConvSec}} \ge h_{\textrm{min}}$}
\If{$w_{\lambda}^{(i+1)} w_{\lambda}^{(i)} > 0$}
\State $Y_{\textrm{cur}}=Y^{(i+1)}$
\Statex 
\Else
\State Continuation of secondary branch terminated.
\State ConvSec=false
\State $Y_{\textrm{cur}}=Y^{(i)}$
\State $Y^{\textrm{final}}=Y^{(i)}$
\EndIf
\Else
\State Continuation of secondary branch terminated.
\Statex
\State ConvSec=false
\Statex
\State $Y_{\textrm{cur}}=Y^{(i)}$
\Statex
\State $Y^{\textrm{final}}=Y^{(i)}$
\EndIf
\EndWhile
\EndIf
\Statex
\If{$||Y_u^{(i+1)}-x_u^{(i+1)}|| \ge 2 \delta_{\textrm{crit}}$}
\State Empty secondary branch
\Statex
\State BREAK and continue the standard Moore-Penrose continuation  
\State algorithm from point $x^{(i+1)}$ (with tangent vector $v^{(i+1)}$). In this case,  
\State we have passed over a bifurcation point and decide to keep following 
\State the principal branch. Figure~\ref{fig:HorizTurnMethod_BifPt} illustrates this scenario.
\EndIf

\EndWhile
\Statex
\Statex Step 3 : If the secondary branch is not empty, flip the points on this branch and add it to the main branch. Then continue the standard Moore-Penrose continuation algorithm from point $Y^{(i)}$ with tangent vector $-w^{(i)}$ (the tangent vector is of opposite sign so that the algorithm keeps pursuing the solution curve in the right direction).

\end{algorithmic}
\end{algorithm}

\end{itemize}

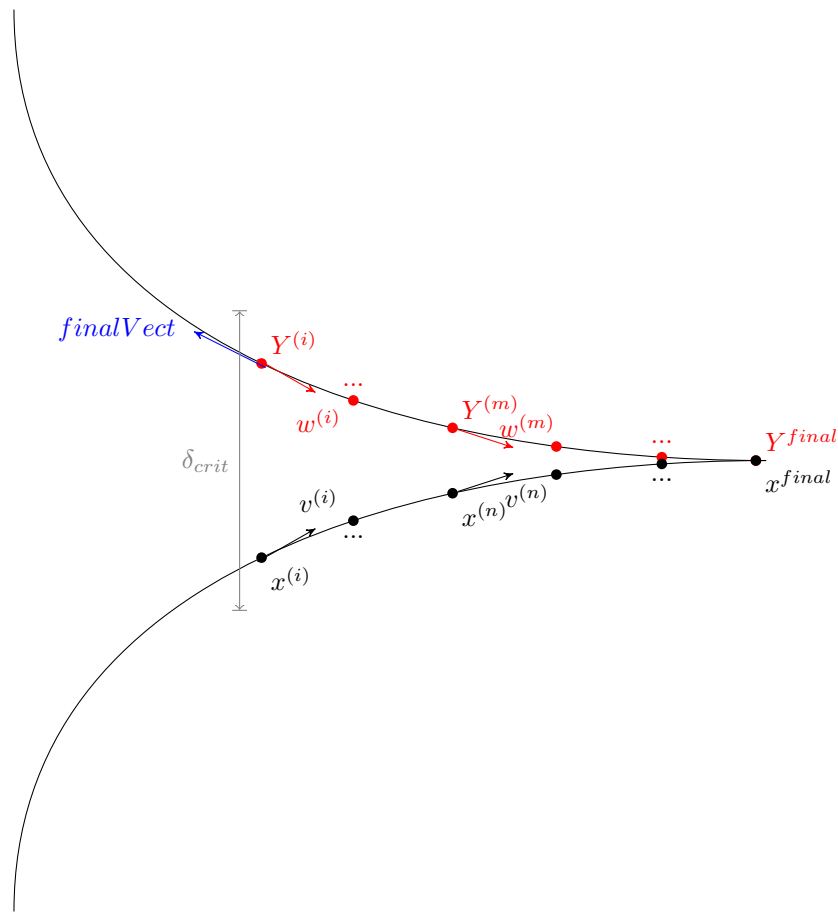
\begin{figure}[!htbp]
    \centering

\begin{tikzpicture}[scale = 2]
\draw (0,3) to[out = 270, in = 180]
node[anchor=south west,pos=0.5, red] {$Y^{(i)}$} 
node[fill,red,circle,inner sep=0pt,minimum size=4pt,pos = 0.5](point2){}
node[pos=0.55] (vect2) {}
node[pos=0.43] (finalvect) {}
node[fill,red,circle,inner sep=0pt,minimum size=4pt,pos = 0.6]{}
node[anchor=south,pos=0.6, above = 2pt, red] {$...$} 
node[fill,red,circle,inner sep=0pt,minimum size=4pt,pos = 0.7](point2n){}
node[anchor=south west,pos=0.7, red] {$Y^{(m)}$} 
node[pos=0.75] (vect2n) {}
node[fill,red,circle,inner sep=0pt,minimum size=4pt,pos = 0.8]{}
node[fill,red,circle,inner sep=0pt,minimum size=4pt,pos = 0.9]{}
node[anchor=south,pos=0.9, above = 2pt, red] {$...$} 
node[fill,red,circle,inner sep=0pt,minimum size=4pt,pos = 0.99](point2f){}
node[anchor=south west,pos=0.99, red] {$Y^{final}$} 
(5,0);

\draw (0,-3) to[out = 90, in = 180] 
node[anchor=north west,pos=0.5] {$x^{(i)}$} 
node[fill,circle,inner sep=0pt,minimum size=4pt,pos = 0.5](point){}
node[pos=0.55] (vect) {}
node[fill,circle,inner sep=0pt,minimum size=4pt,pos = 0.6]{}
node[anchor=north,pos=0.6, below = 2pt] {$...$} 
node[fill,circle,inner sep=0pt,minimum size=4pt,pos = 0.7](pointn){}
node[anchor=north west,pos=0.7] {$x^{(n)}$} 
node[pos=0.75] (vectn) {}
node[fill,circle,inner sep=0pt,minimum size=4pt,pos = 0.8]{}
node[fill,circle,inner sep=0pt,minimum size=4pt,pos = 0.9]{}
node[anchor=north,pos=0.9, below = 2pt] {$...$} 
node[anchor=north west,pos=0.99] {$x^{final}$} 
node[fill,circle,inner sep=0pt,minimum size=4pt,pos = 0.99](pointf){}
(5,0);

\path (point2.north west)
           edge[-stealth',red] node[below = 2pt,pos=1.05] {${{w}^{(i)}}$} (vect2.south east);  

\path (point2n)
           edge[-stealth',red] node[above = 2pt,pos=1.25] {${{w}^{(m)}}$} (vect2n.south east);  

\path (point2.south east)
           edge[-stealth',blue] node[left = 2pt,pos=1.05] {${finalVect}$} (finalvect.west);  
           
\path (point.south west)
           edge[-stealth',black] node[above = 2pt,pos=1.05] {${{v}^{(i)}}$} (vect.north east);  

\path (pointn)
           edge[-stealth',black] node[below = 2pt,pos=1.25] {${{v}^{(n)}}$} (vectn.north east); 

\draw[|<->|,gray] (1.5,-1) -- (1.5,1)
node [anchor=east,pos=0.5] {$\delta_{crit}$};
\end{tikzpicture}
\caption{Illustration of the horizontal turning point method}
\label{fig:HorizTurnMethod} 
\end{figure}

\begin{figure}[!htbp]
    \centering
\begin{tikzpicture}[scale = 1]
\draw (-5,-5) -- (5,5)
node[anchor=north west,pos=0.31] {$x^{(i)}$} 
node[fill,circle,inner sep=0pt,minimum size=4pt,pos = 0.31](point){}
node[pos=0.39] (vect) {}
node[fill,circle,inner sep=0pt,minimum size=4pt,pos = 0.41]{}
node[fill,circle,inner sep=0pt,minimum size=4pt,pos = 0.51]{}
node[fill,circle,inner sep=0pt,minimum size=4pt,pos = 0.61]{}
node[sloped,anchor=north,pos=0.61, above = 2pt] {$...$} 
node[fill,circle,inner sep=0pt,minimum size=4pt,pos = 0.71]{}
node[fill,circle,inner sep=0pt,minimum size=4pt,pos = 0.81]{}
node[fill,circle,inner sep=0pt,minimum size=4pt,pos = 0.9](pointf){}
node[anchor=west,pos=0.9] {$x^{final}$} 
node[pos=0.98] (finalvect) {}
;
\draw (-5,5) -- (5,-5)
node[red,anchor=south west,pos=0.31] {$Y^{(i)}$} 
node[red,fill,circle,inner sep=0pt,minimum size=4pt,pos = 0.31](point2){}
node[red,pos=0.39] (vect2) {}
node[red,fill,circle,inner sep=0pt,minimum size=4pt,pos = 0.41]{}
node[red,fill,circle,inner sep=0pt,minimum size=4pt,pos = 0.51]{}
node[red,fill,circle,inner sep=0pt,minimum size=4pt,pos = 0.61]{}
node[red,sloped,anchor=south,pos=0.61, above = 2pt] {$...$} 
node[red,fill,circle,inner sep=0pt,minimum size=4pt,pos = 0.71]{}
node[red,fill,circle,inner sep=0pt,minimum size=4pt,pos = 0.81]{}
node[red,fill,circle,inner sep=0pt,minimum size=4pt,pos = 0.9](point2f){}
node[red,anchor=south west,pos=0.9] {$Y^{final}$} 
;

\path (point)
           edge[-stealth',black] node[above = 2pt,pos=1.05] {${{v}^{(i)}}$} (vect);  
\path (pointf)
           edge[-stealth',blue] node[left = 2pt,pos=1.05] {${finalVect}$} (finalvect);  
\path (point2)
           edge[-stealth',red] node[anchor = south,right = 2pt] {${{w}^{(i)}}$} (vect2);  

\draw[|<->|,gray] (-1.95,-2) -- (-1.95,2)
node [anchor=east,pos=0.5] {$\delta_{crit}$};
\draw[|<->|,gray] (3.95,-4) -- (3.95,4)
node [anchor=west,pos=0.5] {$2\times\delta_{crit}$};
\end{tikzpicture}
\caption{Illustration of the horizontal turning point method in the case of a bifurcation point}
\label{fig:HorizTurnMethod_BifPt} 
\end{figure}
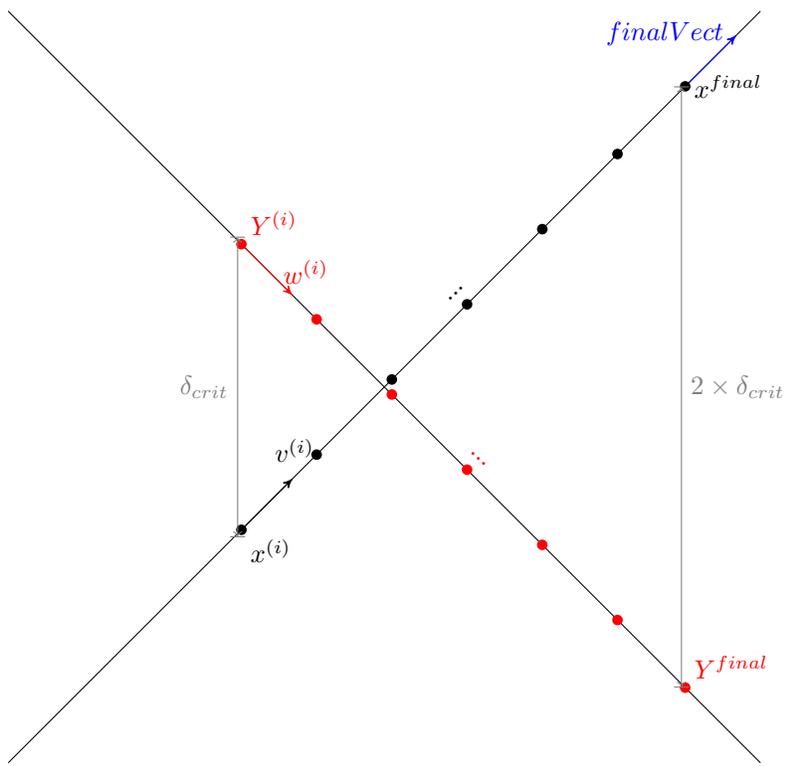

\section{VALIDATION OF THE PROPOSED APPROACH}\label{sec:validation}
To show that our proposed strategy works well, let us consider again the numerical examples of Section~\ref{sec:numex}.

\begin{figure}[!htbp]
\begin{center}
  \begin{tabular}{cc}
   \includegraphics[height=5cm]{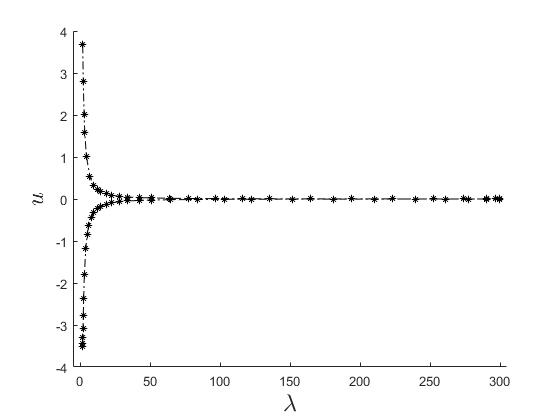} &
   \includegraphics[height=5cm]{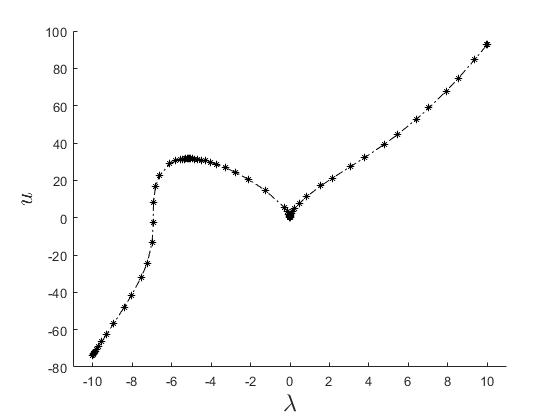} \\
   a) Function $F_a$ & b) Function $F_b$ 
   
  \end{tabular}
  \end{center} 
\caption{Numerical results for functions $F_a$ and $F_b$ using our proposed strategy}
\label{fig:SolutionAlgo_fonctionFaFb}           
\end{figure}

\begin{figure}[!htbp]
\begin{center}
   \includegraphics[height=6cm]{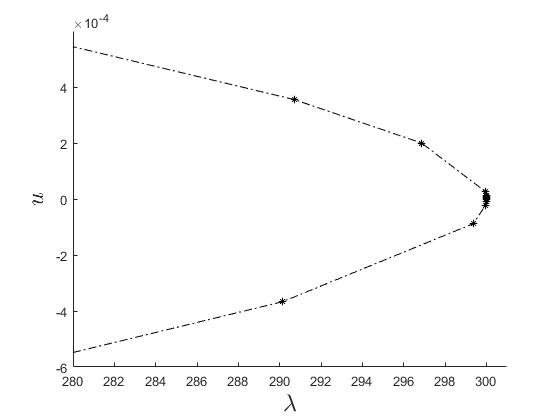}
    \end{center} 
\caption{Numerical results for function $F_a$ using our proposed strategy : zoom in the region of the horizontal limit point}
\label{fig:SolutionAlgoZoom_fonctionFa}           
\end{figure}

Using the standard Moore-Penrose continuation algorithm, both examples of Section~\ref{subsec:case1} had failed to converge past the critical points. Figure~\ref{fig:SolutionAlgo_fonctionFaFb} now shows the results using our proposed approach. In both cases, $c_{\textrm{min}}$ was chosen as $0.95$, $\epsilon_{\lambda}$ was equal to $1 \times 10^{-5}$, $\epsilon_{\lambda}^{\ast}$ was equal to 0.2, $\epsilon_{\textrm{diff}}$ was equal to $1 \times 10^{-7}$ and the deflation technique was applied every $N=5$ iterations. For $F_a$, we used $\delta_{\textrm{maxL}}=30$, $\delta_{\textrm{maxU}}=1.6$ and $\delta_{\textrm{crit}}=2$, while for $F_b$, the values $\delta_{\textrm{maxL}}=1$, $\delta_{\textrm{maxU}}=12$ and $\delta_{\textrm{crit}}=15$ were chosen. As shown, both simulations were completed successfully. To better illustrate the numerical results for function $F_a$, Figure~\ref{fig:SolutionAlgoZoom_fonctionFa} shows a zoom in the critical region. As can be seen, the algorithm calculates many points near the horizontal limit point, which is the more difficult region to compute, and then less points elsewhere.

\begin{figure}[!htbp]
\begin{center}
  \begin{tabular}{cc}
   \includegraphics[height=5cm]{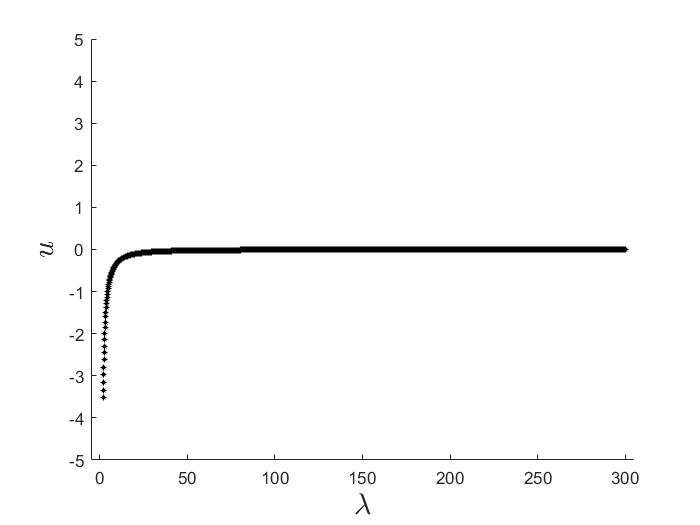} &
   \includegraphics[height=5cm]{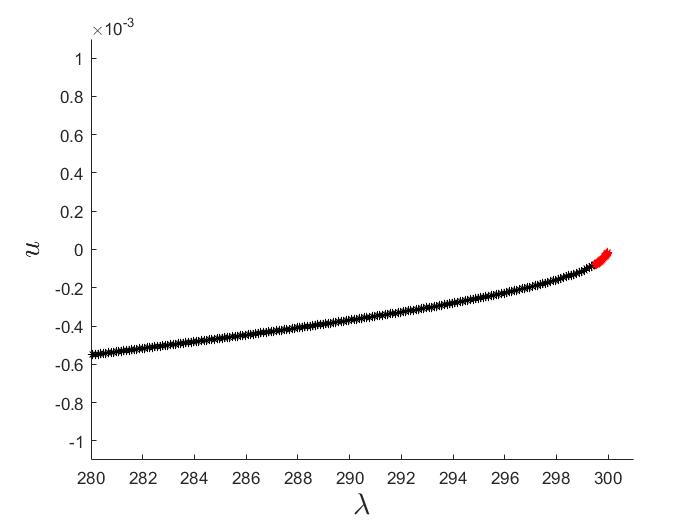} \\
   a) Overall view & b) Zoom in the critical region
   
  \end{tabular}
  \end{center} 
\caption{Numerical results for functions $F_a$ using the Moore-Penrose algorithm with added angle control}
\label{fig:SolutionAC_fonctionFa}           
\end{figure}

To compare our approach with the approach presented in Ligursk\'{y} and Renard~\cite{LigRen2014}, which includes angle control and a simple tangent switch algorithm, we have implemented their approach for function $F_a$. As can be seen in Figure~\ref{fig:SolutionAC_fonctionFa}, the algorithm is not able to complete the simulation adequately, but instead begins to backtrack on the initial region of the solution curve (shown by the red points). By using the same values of parameters as was proposed in their paper, we can also see that adding angle control during the entire simulation with a very strict value for $c_{\textrm{min}}$ is extremely costly as it leads to the calculation of a great number of points. The advantage of using the deflation technique, as well as the other key ingredients presented in Section~\ref{sec:keying}, is that the value of $c_{\textrm{min}}$ in our algorithm does not need to be as strict, which reduces significantly the number of points that will be calculated on the solution curve.

\begin{figure}[!htbp]
\begin{center}
  \begin{tabular}{cc}
   \includegraphics[height=5cm]{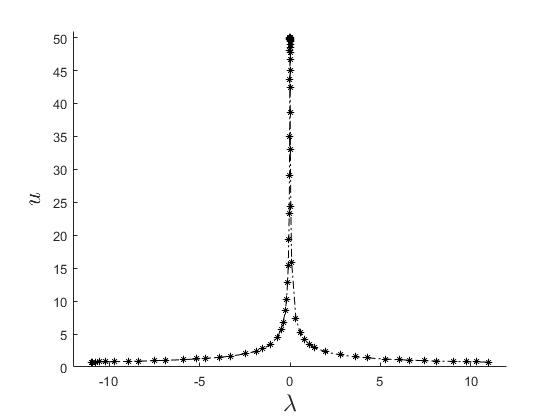} &
   \includegraphics[height=5cm]{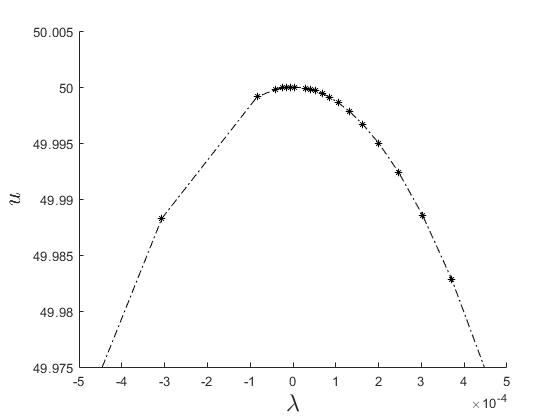} \\
   a) Overall view & b) Zoom in the critical region 
   
  \end{tabular}
  \end{center} 
\caption{Numerical results obtained in the case of $F_c$ using our proposed strategy}
\label{fig:SolutionAlgo_fonctionFc}           
\end{figure}

Let us now consider the example of Section~\ref{subsec:case2} for which the standard Moore-Penrose continuation algorithm started to backtrack on the solution curve. The results using our proposed strategy, with the same parameter values as before for $c_{\textrm{min}}$, $\epsilon_{\lambda}$, $\epsilon_{\lambda}^{\ast}$, $\epsilon_{\textrm{diff}}$ and $N$, but with $\delta_{\textrm{maxL}}=1$, $\delta_{\textrm{maxU}}=10$ and $\delta_{\textrm{crit}}=12.5$, are illustrated in Figure~\ref{fig:SolutionAlgo_fonctionFc}. As can be seen, our algorithm was able to complete the simulation successfully.

\begin{figure}[!htbp]
\begin{center}
  \begin{tabular}{c}
   \includegraphics[height=5cm]{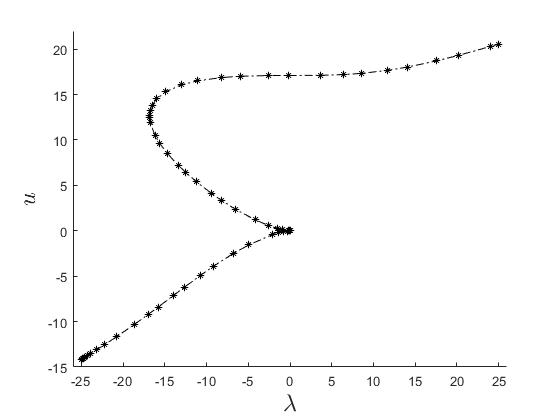} 
     \end{tabular}
  \end{center} 
\caption{Numerical results obtained in the case of $F_d$ using our proposed strategy}
\label{fig:SolutionAlgo_fonctionFd}       
\end{figure}

As for the example in Section~\ref{subsec:case3}, where the standard Moore-Penrose continuation algorithm neglected the entire difficult region of the solution curve by converging to a point much farther along the path, we have once again used our proposed strategy to see if better results can be achieved. Figure~\ref{fig:SolutionAlgo_fonctionFd} shows the results obtained using $\delta_{\textrm{maxL}}=4$, $\delta_{\textrm{maxU}}=1.6$ and $\delta_{\textrm{crit}}=3$. As can be seen, the algorithm no longer neglects the horizontal cusp and is able to complete the simulation without difficulties.

\begin{figure}[!htbp]
\begin{center}
  \begin{tabular}{c}
   \includegraphics[height=5cm]{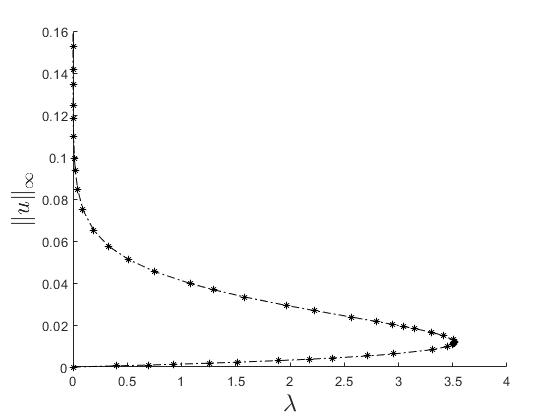} 
     \end{tabular}
  \end{center} 
\caption{Numerical results obtained in the case of the modified Bratu problem with $\gamma=100$ using our proposed strategy}
\label{fig:SolutionAlgo_BratuMod}       
\end{figure}

\begin{figure}[!htbp]
\begin{center}
  \begin{tabular}{c}
   \includegraphics[height=5cm]{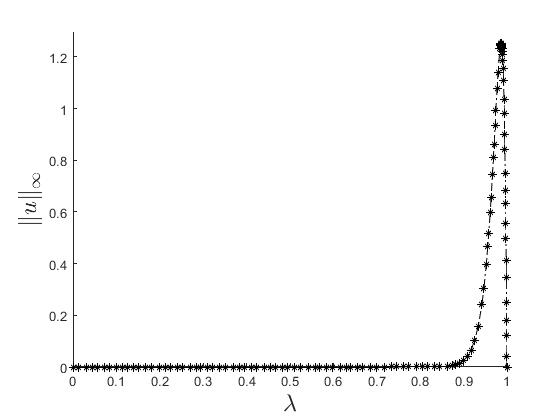} 
     \end{tabular}
  \end{center} 
\caption{Numerical results obtained in the case of the manufactured solution problem using our proposed strategy}
\label{fig:SolutionAlgo_SolMan}       
\end{figure}

Our algorithm was also tested on the finite element problems presented in Section~\ref{subsec:finelemprob}. In both cases, the simulation was completed successfully as can be seen in Figures~\ref{fig:SolutionAlgo_BratuMod} and~\ref{fig:SolutionAlgo_SolMan}. The numerical results were obtained using $\delta_{\textrm{maxL}}=0.1$, $\delta_{\textrm{maxU}}=0.02$ and $\delta_{\textrm{crit}}=0.025$ for the modified Bratu problem and $\delta_{\textrm{MaxL}}=0.02$, $\delta_{\textrm{MaxU}}=0.2$ and $\delta_{\textrm{crit}}=0.25$ for the manufactured solution problem.

All of the numerical examples of Section~\ref{sec:numex} were therefore completed successfully using our proposed approach. Our algorithm is much more robust than the standard Moore-Penrose continuation method and is able to compute the solution curve of challenging problems without difficulty.

\begin{figure}[!htbp]
\begin{center}
  \begin{tabular}{c}
   \includegraphics[height=5cm]{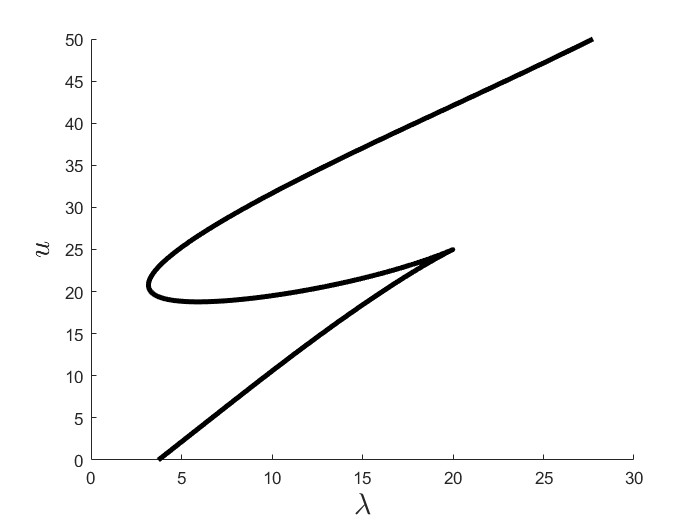} 
     \end{tabular}
  \end{center} 
\caption{Solution curve for function $F_e$}
\label{fig:SolutionExacte_fonctionFe}           
\end{figure}

\begin{figure}[!htbp]
\begin{center}
  \begin{tabular}{cc}
   \includegraphics[height=5cm]{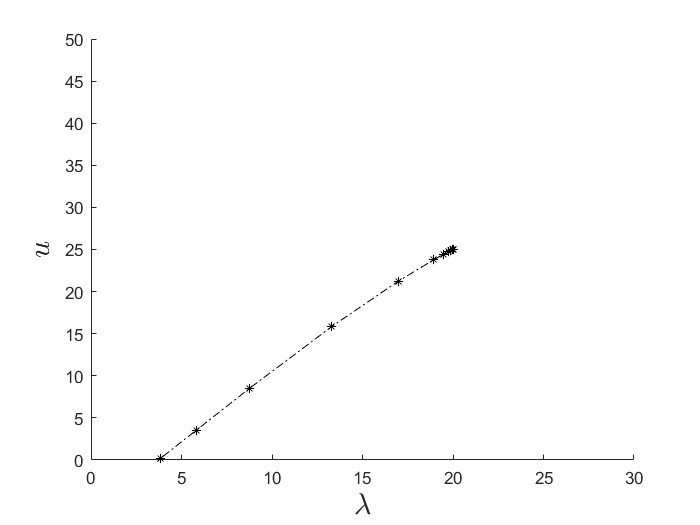} &
   \includegraphics[height=5cm]{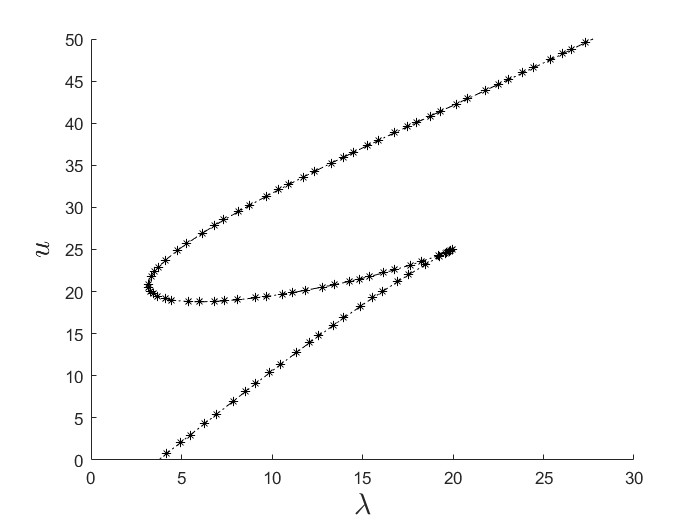} \\
   a) Standard algorithm & b) Proposed approach 
   
  \end{tabular}
  \end{center} 
\caption{Numerical results obtained in the case of $F_e$}
\label{fig:SolutionAlgo_fonctionFe}           
\end{figure}

To test our approach on a problem where the limit state on the solution curve appears at an angle, a behaviour that could be encountered in more complex problems, let us consider the function $F_e :  \mathbb{R}^2 \longrightarrow \mathbb{R}$ defined respectively by : 
$$F_e(u, \lambda) = -500(\lambda-u-5)^2-10(u-20)^3+0.1(\lambda-u-5)^5$$
The solution curve for $F_e(u, \lambda)=0$ is shown in Figure~\ref{fig:SolutionExacte_fonctionFe} and is very similar to the solution curve obtained in the case of shell buckling (see~\cite{Far2006}). As can be seen in Figure~\ref{fig:SolutionAlgo_fonctionFe}, the standard approach stops converging at the critical point while our proposed approach has no difficulty to complete the simulation and trace the entire solution curve.

\begin{figure}[!htbp]
\begin{center}
  \begin{tabular}{c}
   \includegraphics[height=5cm]{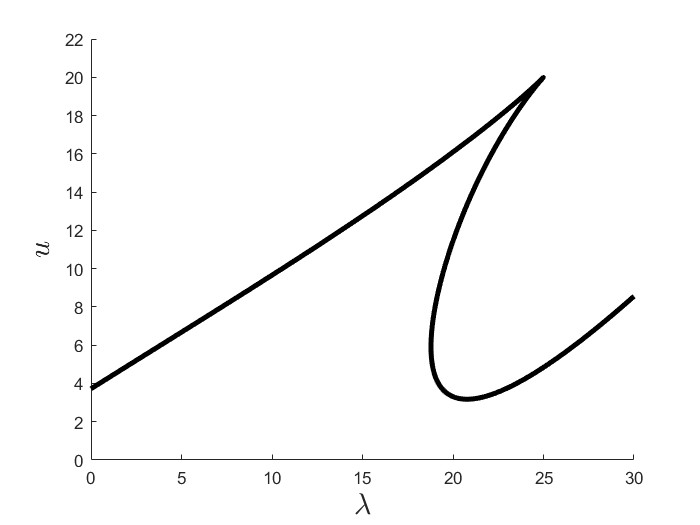} 
     \end{tabular}
  \end{center} 
\caption{Solution curve for the inverse of function $F_e$}
\label{fig:SolutionExacte_fonctionFe_inv}           
\end{figure}

\begin{figure}[!htbp]
\begin{center}
  \begin{tabular}{cc}
   \includegraphics[height=5cm]{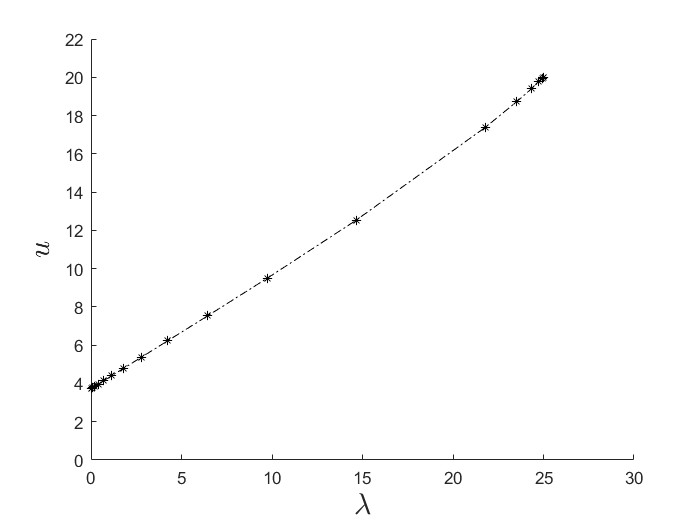} &
   \includegraphics[height=5cm]{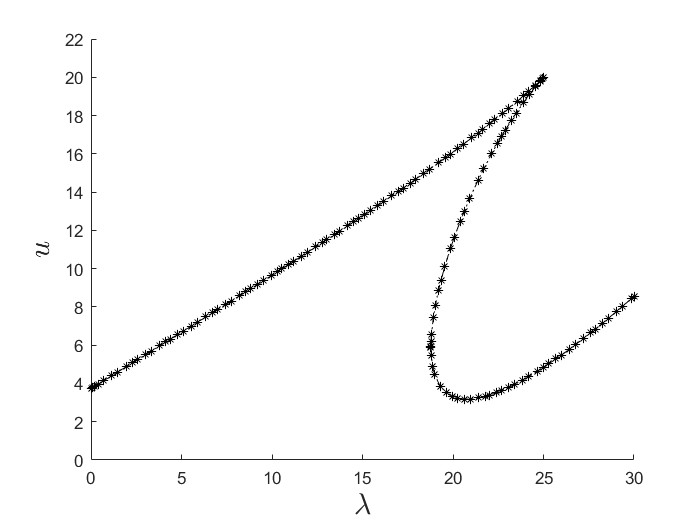} \\
   a) Standard algorithm & b) Proposed approach 
   
  \end{tabular}
  \end{center} 
\caption{Numerical results obtained in the case of the inverse of $F_e$}
\label{fig:SolutionAlgo_fonctionFe_inv}           
\end{figure}

We can also look at what happens if we invert the graph so that the solution curve now exhibits a slightly different scenario in terms of the critical point. The solution curve is shown in Figure~\ref{fig:SolutionExacte_fonctionFe_inv} while Figure~\ref{fig:SolutionAlgo_fonctionFe_inv} shows the numerical results for both the standard algorithm and our proposed approach. Again, our approach leads to very good results while the standard approach is not able to converge past the critical point.

\subsection{General remarks}

General remarks can be made in regards to the proposed algorithm for the Moore-Penrose continuation method.

\begin{itemize}

\item The choice of the parameter values $\delta_{\textrm{maxL}}$, $\delta_{\textrm{maxU}}$ and $\delta_{\textrm{crit}}$ depends on the problem we are solving. We generally choose them in terms of the span of $u$ and $\lambda$. In most of our problems, we chose $\delta_{\textrm{maxL}}$ as the span of $\lambda$ (of the difficult region) divided by 10, $\delta_{\textrm{maxU}}$ as the span of $u$ (of the difficult region) divided by 5 and $\delta_{\textrm{crit}}$ as the span of $u$ (of the difficult region) divided by 4. This can be used as a general guideline for choosing these values. 
Let us note that $\delta_{\textrm{crit}}$ should always be greater than $\delta_{\textrm{maxU}}$.

\item The proposed approach does introduce other new parameters such as  $c_{min}$ and $N$. From our sensitivity analysis, a good range of values can be used for these parameters. As a guideline, one could choose $c_{min}$ in the interval $[0.9, 1[$ and $N$ lower or equal to 10. By using these values, all of our numerical tests were completed successfully. Choosing a lower value for $c_{min}$ in the given range will lead to the calculation of less points in critical regions while choosing a value closer to $1$ will lead to the calculation of more points. The severity of the critical point could help determine what value should be chosen for this parameter. For all our examples, we chose $c_{min}=0.95$ and did not encounter situations where this value needed to be closer to $1$. As for $N$, if it is chosen too large, we might not be able to identify in advance, for example, that an horizontal limit point is approaching. That is why we chose $N=5$ for our numerical problems.

\item During the validation process, we also considered multidimensional problems where the solution curve for at least one of the variables in terms of $\lambda$ exhibited a critical point. The proposed approach led to good results in these examples as well.

\item In the \textit{Vertical turning point method}, Newton's method is used to find $Z^{\ast}$. As the algorithm was able to converge to $x_u^{(i)}$ and $\Delta_{\lambda}$ is chosen small, one can reasonably expect convergence at this step. If however Newton's method does not converge, reducing slightly the value of $\Delta_{\lambda}$ should help with the situation.

\end{itemize}

\section{CONCLUSION}
In this work, we have showed that the standard Moore-Penrose continuation method can lead to undesired results near critical points. In particular, the algorithm can diverge, backtrack on a part of the solution curve that has already been computed or omit important regions of the solution curve. This paper therefore presents a more robust approach for the computation of the solution curve when difficulties occur. The approach is a modified version of the Moore-Penrose continuation method and includes four key ingredients : the deflated continuation algorithm, angle control, a control on the distance between consecutive converged points as well as the monitoring of the sign of $v_{\lambda}$. From the validation tests considered, we have seen that this algorithm leads to significantly better results and is very promising.

\section*{Acknowledgements}
The authors wish to acknowledge the financial support of the Faculty of Science of the Université de Moncton and the Natural Sciences and Engineering Research Council of Canada (NSERC) (Grant number RGPIN-2017-05099).

%\newpage

\bibliographystyle{plain} % Style bibliographique auteur-date. Pour un style avec numéro, remplacer par : \bibliographystyle{plain}
\bibliography{banquebiblio,biblio_afortin} %Pour faire sa bibliographie avec Bibtex

\begin{thebibliography}{10}

\bibitem{AllGeo1990}
E.L. Allgower and K.~Georg.
\newblock {\em An Introduction to Numerical Continuation Methods}, volume~13 of
  {\em Springer series in computational mathematics}.
\newblock Springer-Verlag, New York, 1990.

\bibitem{Bra1914}
G.~Bratu.
\newblock Sur les équations intégrales non linéaires.
\newblock {\em Bulletin de la Société Mathématique de France}, 42:113--142,
  1914.

\bibitem{ChaForFor2010}
{\'E}.~Chamberland, A.~Fortin, and M.~Fortin.
\newblock Comparison of the performance of some finite element discretizations
  for large deformation elasticity problems.
\newblock {\em Computers \& Structures}, 88(11-12):664--673, 2010.

\bibitem{CheSch1990}
Z.~Chen and H.L. Schreyer.
\newblock A numerical solution scheme for softening problems involving total
  strain control.
\newblock {\em Computers \& Structures}, 37(6):1043--1050, 1990.

\bibitem{Cri1981}
M.~A. Crisfield.
\newblock A fast incremental/iterative solution procedure that handles snap
  through.
\newblock {\em Computers \& Structures}, 13:55--62, 1981.

\bibitem{DhoGovKuz2003}
A.~Dhooge, W.~Govaerts, and Y.~A. Kuznetsov.
\newblock {MATCONT} : a {MATLAB} {p}ackage for {n}umerical {b}ifurcation
  {a}nalysis of {ODE}s.
\newblock {\em ACM Transactions on Mathematical Software}, 29(2), 2003.

\bibitem{FarBeeBir2016}
P.~E. Farrell, C.~H.~L. Beentjes, and \'{A}. Birkisson.
\newblock The computation of disconnected bifurcation diagrams.
\newblock {\em arXiv:1603.00809}.

\bibitem{FarBirFun2015}
P.~E. Farrell, \'{A}. Birkisson, and S.~W. Funke.
\newblock Deflation techniques for finding distinct solutions of nonlinear
  partial differential equations.
\newblock {\em SIAM Journal of Scientific Computing}, 37(4):2026 -- 2045, 2015.

\bibitem{Far2006}
M.~Farshad.
\newblock {\em Plastic Pipe Systems: Failure Investigation and Diagnosis}.
\newblock Elsevier Science, 1 edition, 2006.

\bibitem{Fri1984}
I.~Fried.
\newblock Orthogonal trajectory accession to the non-linear equilibrium curve.
\newblock {\em Computer Methods in Applied Mechanics and Engineering}, 47:283
  -- 298, 1984.

\bibitem{Gut2004}
M.A. Gutiérrez.
\newblock Energy release control for numerical simulations of failure in
  quasi-brittle solids.
\newblock {\em International Journal for Numerical Methods in Engineering},
  20(1):19--29, 2004.

\bibitem{JacSch2002}
J.~Jacobsen and K.~Schmitt.
\newblock The {L}iouville-{B}ratu-{G}elfand problem for radial operators.
\newblock {\em Journal of Differential Equations}, 184(1):283--298, 2002.

\bibitem{LigRen2014}
T.~Ligursk\'{y} and Y.~Renard.
\newblock A continuation problem for computing solutions of discretised
  evolution problems with application to plane quasi-static contact problems
  with friction.
\newblock {\em Computer methods in applied mechanics and engineering}, 280:222
  -- 262, 2014.

\bibitem{LegDeiFor2015}
S.~Léger, J.~Deteix, and A.~Fortin.
\newblock A {M}oore-{P}enrose continuation method based on a {S}chur complement
  approach for nonlinear finite element bifurcation problems.
\newblock {\em Computers \& Structures}, 152:173--184, 2015.

\bibitem{LegForTib2014}
S.~Léger, A.~Fortin, C.~Tibirna, and M.~Fortin.
\newblock An updated {L}agrangian method with error estimation and adaptive
  remeshing for very large deformation elasticity problems.
\newblock {\em International Journal for Numerical Methods in Engineering},
  100(13):1006--1030, 2014.

\bibitem{LegPep2016}
S.~Léger and A.~Pepin.
\newblock An updated {L}agrangian method with error estimation and adaptive
  remeshing for very large deformation elasticity problems : the
  three-dimensional case.
\newblock {\em Computer Methods in Applied Mechanics and Engineering}, 309:1 --
  18, 2016.

\bibitem{PohRamBis2014}
T.~Pohl, E.~Ramm, and M.~Bischoff.
\newblock Adaptive path following schemes for problems with softening.
\newblock {\em Finite Elements in Analysis and Design}, 86:11 -- 22, 2014.

\bibitem{Rik1979}
E.~Riks.
\newblock An incremental approach to the solution of snapping and buckling
  problems.
\newblock {\em International Journal of Solids and Structures}, 15:529--551,
  1979.

\bibitem{VanProSim2013}
B.~Vandoren, K.~De~Proft, A.~Simone, and L.J. Sluys.
\newblock A novel constrained {LA}rge {T}ime {IN}crement method for modeling
  quasi brittle failure.
\newblock {\em Computer Methods in Applied Mechanics and Engineering}, 265:148
  -- 162, 2013.

\end{thebibliography}

\end{document}